%% file: les.tex
\documentclass[hein,final,leqno]{shein}
\usepackage{latexsym}
\usepackage{graphics}
\usepackage{epsfig}
\usepackage{epsf}
\usepackage{eucal}
\usepackage{pstcol}
\usepackage{amsfonts}
\usepackage{amstext}
\usepackage{amssymb}
\usepackage{amsmath}
\usepackage{amsbsy}
\usepackage{amsxtra}
\usepackage{amscd}
\shnewtheorem*{udefinition}{Definition}{\bfseries\itshape}{} 
\shnewtheorem*{uremark}{Remark}{\bfseries\itshape}{}   
\shnewtheorem*{uremarks}{Remarks}{\bfseries\itshape}{} 
\shnewtheorem*{utheorem}{Theorem.}{\bfseries}{} 
\renewcommand{\abstractname}{Abstract.}
\begin{document}
\title{Large eddy approximation of turbulent flow in DSC schemes}
\author{Steffen Hein
}                     
%
%
\institute{Steffen Hein 
\hfill DE-83043 Bad Aibling,
Germany\hfill}
%
%
%
\titlerunning{Large eddy approximation in DSC schemes}
\authorrunning{Steffen Hein}

\maketitle 
%
\vspace{-12pt}
\input{abs.tex}   
%
\markboth{{\normalsize \textsc{Steffen Hein}}}
{{\normalsize \textsc{Large eddy approximation}}}

\normalsize
\vspace{-6pt}
\input{qtn.tex}   
\vspace{-19pt}
\input{int.tex}   
\vspace{-8pt}
\input{sec1.tex}  
\vspace{-8pt}
\input{sec2.tex}  
\vspace{-8pt}
\input{sec3.tex}  
\vspace{-8pt}
\input{sec4.tex}  
\vspace{-8pt}
\input{sec5.tex}  
\vspace{-8pt}
\input{ccl.tex}   
\vspace{-8pt}

\input{ref.tex}
\hrulefill
\newline
\textsc{Steffen Hein};\;
DE-83043 Bad Aibling, Germany
\newline E-mail address:\; steffen.hein@bnro.de
\end{document}

%% file: abs.tex
\begin{abstract}
Large eddy simulation of turbulent flow is given a natural setting
within the DSC framework of computational fluid dynamics. Periodic cellular
coarse-graining prevents the nodal flow from piling up
and preserves its large patterns. The coarsening operations are
consistent with the near-field interaction principle of DSC and
-~therefore~- uncomplicated at boundaries.
Numerical examples validate the approach.
{\vspace{2pt}\hfill
\textbf{MSC-classes}:\;\textnormal{65C20,\,65M06,\,76D05}}

\medskip
\noindent
{\bf Keywords: \small{Navier-Stokes equations, turbulence models, 
large eddy \linebreak simulation, DSC schemes}}
\medskip
\vspace{-10pt}
\end{abstract}

%% file: qtn.tex
\begin{quote}\notag
\small
{
\makebox{
\textnormal{L'acte de la nature est simple, et la complexit\'e
qu'elle para\^{\i}t avoir}}
\newline
\makebox{
\textnormal{construite ... n'est que l'entrecroisement sans fin
des antagonismes}}
\newline
\makebox{
\textnormal{qui se sont neutralis\'es}}
\vspace{1pt}
\begin{flushright}
\textit{Henri Bergson}
\end{flushright}
}
\end{quote}

%% file: int.tex
\section{Introduction}\label{S:int}
Discrete schemes have to master the cursed situation that they
first create by breaking space into pieces.
Mesh cell systems, however fine, can never be perfect substitute
for continuous space. Already in linear algorithms they produce 
artefacts, such as unphysical \emph{spurious solutions},
and more harmful things happen in the non-linear case.

In fluid dynamics, when the non-linear spectral transfer properties
of the Navier-Stokes equations
(the \emph{energy cascade}; cf. \textsc{Pope} \cite{Po})
come into play beyond the transition to turbulence,
eddies are locally excited, down to very short scale;
cf. \textsc{Kolmogorov} \cite{Ko}.
Eddies smaller than the cell size cannot be properly resolved by a mesh
of realistic coarseness and thus tend to induce local fluctuations that
artificially pile up. Large eddy simulation (LES) aims to discard
such divergences through regularizing the flow by means of a suitable
averaging filter, and to screen its essential \emph{large patterns}
in this way.

LES stands for the most promising line in modern turbulence modelling.
The name has been coined by \textsc{Deardorff} in 1970,
who first applied LES methods to turbulent channel flow \cite{De}.
Important elements had yet been kept ready, formerly.
Already \textsc{Reynolds} \cite{Re} used temporal flow averages
that as \emph{Reynolds averages}
play still a r\^ole in conventional turbulence models.
In 1922, \textsc{Richardson} \cite{Ri} proposed mesh cell averages
for smoothing down local fluctuations. Such so called \emph{box filters}
stand today for one line of LES. Spatial filters of varied type indeed
classify modern LES methods. Excellent insight into the state of the
art is gained from the recent book of \textsc{Berselli} et al. \cite{BIL},
which not only exposes the elaborate mathematical framework but also
addresses many items of conventional turbulence modelling 
(the prominent $k\!-\epsilon$ model of
\textsc{Launder} and \textsc{Spalding} \cite{LS}, e.g.).

The present paper differs from most related work in the field in that
it sharply separates the natural and technical aspects
of turbulent instability and strictly focusses on the
latter within the \emph{computational context} in hand.
The natural part, viz. essentially (modulo discretization)
the Navier-Stokes equations for viscous incompressible flow,
are thereby taken as a \emph{physical given}.
This is in obvious contrast to widespread use:
In many conventional turbulence models, as in some lines of LES,
turbulent pile-up is controlled by modifying the underlying dynamic
equations or even the physical constants.
Such models use, for instance, Reynolds number dependent
global parameters -~usually an \emph{eddy viscosity}, e.g.~-
or time and space filtered Navier-Stokes equations or any 
(e.g. stochastic) variant.
All this is allowed, of course, as long as it provides numerical
data in harmony with observation. Also, a profound change of the
inner dynamics can very efficiently delimit turbulence and thus
lead to a convergent scheme.
It is yet, obviously, tuning the underlying \emph{physical phenomenon}
rather than eliminating the technical cause of divergences
that are manifestly an \emph{artefact of discretization}\/.
(Turbulence is still enigmatic in many respects, but \emph{divergence}
is certainly not a feature of turbulent flow -~observed in the continuum
of nature, e.g.)

Therefore, the central point of attack in discarding turbulent pile-up
will be in this paper the weaknesses of discretization,
while the Navier-Stokes equations are directly discretized as taken
from physics -~which not only combats the desaster at the origin,
but also yields validated results.

%% file: sec1.tex
\section{DSC schemes}\label{S:sec1}
A mesh cell system is not just poor surrogate continuous space,
but it always introduces its proper structure. 
Cellular meshes, with that we are dealing here,
artificially impose local \emph{cell-boundary duality} 
upon space~\cite{He1}
-~and it is this simple fact that dual scattering channel
(DSC) schemes match in a canonical fashion.

DSC schemes are finite volume methods of a very specific type. 
In technical detail, they are extensively dealt with in~\cite{He1},
so it is sufficient here to recapitulate their essential features
in a crash course like manner.
The algorithm is characterized by a two-step cycle of iteration
which alter\-nately updates the computed fields within cells and on
their interfaces. If the updating instructions are explicit, then
a \emph{near-field interaction} principle gives rise to a scattering
process interpretation known as \emph{Johns' cycle}~\cite{He1}. 
A pair of vectors that represent the same field within a cell and
on its surface essentially constitutes a \emph{scattering channel}.
Equivalently, scattering channels are sometimes defined as pairs of
distributions that 'measure' the field within a cell and on one of
its faces. 
The well known primal DSC scheme is \textsc{Johns}' TLM algorithm,
wherein scattering channels are visualised as transmission
lines~\cite{JoB}.
In many applications, such as computational fluid dynamics, the
TLM picture of wave propagation fails, but one can still go back
to the far more general DSC setup.

Within any \emph{physical interpretation\,}, cf.~\cite{He1}\,,
each scattering channel represents a pair of scalar or vector valued
distributions
$z^p\,=\,(\,p\,,\,Z\,)\,$ and $z^n\,=\,(\,p\,\sptilde,\,Z\,)\,$
which evaluate a physical field $Z\/$ on the surface and in the interior,
respectively, of a mesh cell. A distribution
$\,p\,$ pertinent to a cell face is called a \emph{port} 
and $\,p\,\sptilde\/$ its \emph{nodal image},
and the two are related by pull-back, viz.
\vspace{-8pt}
\begin{equation}\centering\label{1.1}
(\,p\,\sptilde,\,Z\,)\;
=\;(\,p\,\circ\,\sigma\,,\,Z\,)\;
=\;(\,p\,,\,Z\,\circ\,\sigma^{-1}\,)\;,
\end{equation}
for every $\,Z\,$ (\,of class $\,C^{\,\infty}\,$, e.g.\,)\;,
where $\,\sigma\,$ denotes the spatial translation
${\sigma:\mathbb{R}^3\,\to\,\mathbb{R}^3}$ 
that shifts the geometric \emph{node} (i.e. the cell centre) in the
(centre of) the respective face.

It follows that there exists a \emph{scattering channel representation\,} 
of DSC states:
If $M$ denotes the mesh cell system, $\,{\partial\zeta}\,$
the boundary of cell $\,{\zeta\in M}\,$
($\,{\partial\zeta}\,$ is naturally identified with the set of ports
with nodal image in $\,\zeta\,$), and $L_{\,\zeta}^{\,p}$ the 
span of the $z_{\,\zeta}^{\,p}\;$ (i.e. the linear \emph{range} space),
then each state in the mesh permits a unique representation in the space
\begin{equation}\centering\label{1.2}
P\;:\,=\;\prod\nolimits_{\,\zeta\in M}\;
\prod\nolimits_{\,p\in\partial\zeta}\;
L_{\,\zeta}^{\,p}\,\times\,L_{\,\zeta}^{p\sptilde}\;,
\end{equation}
with canonical projections
${\pi_{\,\zeta}^{\,p,\,n}\,:\,P\,\to\,P_{\,\zeta}^{\,p,\,n}}$
into the port and node components of cell $\zeta\;$.
Furthermore, there is a natural involutary isomorphism 
${nb\,:\,P\,\to\,P\,}$, 
\vspace{-8pt}
\begin{equation}\centering
nb\,:\,(\,z^{\,p}\,,\,z^{\,p\sptilde}\,)\,\mapsto\,
(\,z^{\,p\sptilde},\,z^{\,p}\,)\;,
\end{equation}
called the \emph{node-boundary map}, which hence maps
${P^{\,p}}$ onto ${P^{\,n}}$ and vice versa.
(The cell index $\zeta$ is omitted here and in the following
without danger of confusion.)

It stands to reason that the nodal images of different ports
(e.g. on different faces of the same cell)
can represent the same physical field in a node 
and may even coincide there as distributions -~just as two ports pertinent
to neighbouring cells obviously represent the same field on a common face, 
if they \emph{connect} two channels on it.
In fact, the scattering channel representation of DSC states in \eqref{1.2}
is in general highly \emph{redundant}
-~which may be utilized for process parallelization~\cite{He1}.

Within the algorithm, the port and node components are updated at even
and odd integer multiples, respectively, of half a timestep ~${\tau}\,$ 
and are usually constantly continued as step functions over the
subsequent time intervals of length $\tau\,$.

Note that \emph{existence}, not necessarily explicit \emph{construction}
or \emph{application}, of a scattering channel representation
characterizes DSC schemes, and that less redundant representations
are actually used in many implementations.
In fact, the scattering channel representation is basically a
theoretical means for describing (and deriving) certain DSC properties 
-~such as the following.

A fundamental principle -~closely related to the \textsc{Courant-Levi} 
stability criterion~- is \emph{near-field interaction}.  
It requires that every updated state of a node or face depends 
only states (along with their history) in scattering channels
\emph{connected} to the respective node or face -~the latter 
\emph{here} being for once identified with its adjacent face,
if such exists in any neighbouring cell.

As a consequence of near-field interaction, every DSC process allows
for an interpretation as a \emph{multiple scattering process} in the
following sense. 

Let for any process $\,z\,=\,{\,(\,z^{\,p}\,,\,z^{\,n}\,)(\,t\,)\,}$
\emph{incident\/} and \emph{outgoing fields}
$\,{z_{\,in}^{\,p}}\,$ and $\,{z_{\,out}^{\,n}}\,$ 
be recursively defined as processes in
${\,P^{\,p}\,}$ and ${\,P^{\,n}\,}$, respectively, by setting 
$\,{z_{\,in}^{\,p}\,(\,t\,)\,
:\,=\,z_{\,out}^{\,n}\,(\,t-\frac{\tau}{2}\,)\,
:\,=\,0\;}\,$,\;\;for $\,t\,<\, 0\,$, \newline 
and for $\,0\,\leq\,t\,=\,m\tau\,$;\,$\,m\in\mathbb{N}\,$\::
\vspace{-8pt}
\begin{equation}\centering\label{1.4}
\begin{split}
\begin{aligned}
z_{\,in}^{\,p}\,(\,t\,)\;
&:\,=\; z^{\,p}\,(\,t\,)\,
-\,nb\circ z_{\,out}^{\,n}\,(\,t\,-\,\frac{\tau}{2}\,)\;,\\
z_{\,out}^{\,n}\,(\,t\,+\,\frac{\tau}{2}\,)\;
&:\,=\; z^{\,n}\,(\,t+\,\frac{\tau}{2}\,)\,-\,nb\circ z_{\,in}^{\,p}\,(\,t\,)\;.
\end{aligned}
\end{split}
\end{equation}
Then, at every instant holds \,$\,{z^{\,p}\,(\,t\,)}\,
=\,{nb\circ z_{\,out}^{\,n}\,(\,t\,-\,\frac{\tau}{2}\,)}\,
+\,{z_{\,in}^{\,p}\,(\,t\,)}\,$\,
and \\
$\,{z^{\,n}\,(\,t\,+\,\frac{\tau}{2}\,)}\,
=\,{nb\circ z_{\,in}^{\,p}\,(\,t\,)}\,
+\,{z_{\,out}^{\,n}\,(\,t\,+\,\frac{\tau}{2}\,)}\,$\,.
Also, near-field interaction implies that every state is only a function of
states incident (up to present time $\,t\,$) on scattering channels
connected to the respective node or face.
\newline
More precisely, by induction holds:
\begin{utheorem}
There exists a pair of functions
$\,\mathcal{R}\,$ and $\,\mathcal{C}\,$, 
defined on back in time running sequences of incident and outgoing
fields, respectively,
such that for every cell $\,{\zeta\in M}\,$
the process $\,z_{\,\zeta}^{\,n}\,=\,{\pi_{\,\zeta}^{\,n}\circ z}\,$
complies with
\vspace{-8pt}
\begin{equation}\centering\label{1.5}
z_{\,\zeta}^{\,n}\,(\,t\,+\,\frac{\tau}{2}\,)\;
=\;\mathcal{R}\,(\,(\,z_{\,in}^{\,p}\,(\,t\,-\,\mu\tau\,)\,)
_{\,p\in\partial\zeta\,;\;\mu\in\mathbb{N}}\,)\;
\vspace{-8pt}
\end{equation}
and the port process
$\,z_{\,\zeta}^{\,p}\,=\,\pi_{\,\zeta}^{\,p}\,\circ\,z\,$
satisfies
\vspace{-10pt}
\begin{equation}\centering\label{1.6}
\begin{aligned}
z_{\,\zeta}^{\,p}\,(\,t\,+\,\tau\,)\;
=\;\mathcal{C}\,(\,(\,z_{\,out}^{\,n}\,(\,t\,+\,\frac{\tau}{2}-
\mu\tau\,)\,)_{\,n\,\mid\,\partial\zeta\,;\;\mu\in\mathbb{N}}\,)\;.\\
\hfill\scriptsize{(\;\text{'$\;\mid\,$' short-hand for
'in any of the (\emph{\.1 or 2\,}) cells adjacent to'}\;\;)} 
\end{aligned}
\end{equation}
\end{utheorem}
\begin{uremarks}
\vspace{-8pt}
\hspace{24pt}
\begin{itemize}
\item[(i)]
The statements immediately imply that
$\,z_{\zeta,\,out}^{\,n}\,$ and $\,z_{\zeta,\,in}^{\,p}\,$
are themselves functions of states incident
on connected scattering channels, since 
\vspace{-8pt}
\end{itemize}
\begin{equation}\centering\notag
\begin{aligned}
\quad\; z_{\zeta,\,out}^{\,n}\,(\,t\,+\,\frac{\tau}{2}\,)\,
&=\,\mathcal{R}\,(\,(\,z_{\,in}^{\,p}\,(\,t\,-\mu\tau\,)\,)_{\,p\,\in\,
\partial\zeta\,;\;\mu\in\mathbb{N}}\,)\,
-\,nb\circ\,z_{\zeta,\,in}^{\,p}\,(\,t\,)\quad\text{and}\\
z_{\zeta,\,in}^{\,p}\,(\,t\,)\,
&=\;\mathcal{C}\,(\,(\,z_{\,out}^{\,n}\,(\,t\,-\,\frac{\tau}{2}-\mu\tau\,)\,
)_{\,n\,\mid\,p\,;\;\mu\in\mathbb{N}}\,)\;\;+\\[-3pt]
&\qquad\qquad\qquad\qquad\qquad\qquad
-\;nb\circ\,z_{\zeta,\,out}^{\,p\sptilde}\,(\,t\,-\,\frac{\tau}{2}\,)
\end{aligned}
\vspace{-4pt}
\end{equation}
\begin{itemize}
\item[(ii)]
$\,\mathcal{R}\,$ and $\,\mathcal{C}\,$ are called, respectively, the
\emph{reflection} and \emph{connection\/} maps of the DSC algorithm, and
a \emph{field excitation} may be implicit in $\,\mathcal{C}\,$, cf.~\cite{He1}.
\end{itemize}
\begin{itemize}
\item[(iii)]
Near field interaction implies computational stability, if the reflection
and connection maps are 
\nolinebreak{$\,\alpha$-\emph{passive}}\,, 
i.e. contractive in this sense~\cite{He2}.
\end{itemize}
\end{uremarks}

%% file: sec2.tex
\section{Coarse-graining}\label{S:sec2}
Mesh cell systems import artificial defects into space
-~which yet sometimes carry their remedy in themselves, at least in part: 
For instance, the cellular mesh, by separating the virtually unresolved
cell interior from the coarse-grained cell-boundary \emph{skeleton grid\/},
necessarily fixes -~and hence \emph{delimits\/}~- 
the scale of local resolution. 
A clever design always takes into account the coarsening effect of the
cellular mesh -~which, on the other hand, provides quasi free of charge
a basis for large eddy approximation in the DSC setup of computational
fluid dynamics. 
In fact, already intuitively DSC schemes go along with kind of natural
large pattern approximation that -~in a sense now made precise~-
associates \emph{cell averages} to some nodal quantities.

The following definition is independent of any particular application
and therefore given without reference to fluid dynamics.

\begin{udefinition} (coarse-graining)
\newline
For any port 
${\,p\in\partial\zeta\,}$ with nodal image $\,p\sptilde\,$ 
let ${\,B^{p}\,}$ denote the set
of all ports on $\,\partial\zeta\,$, the nodal images of which
coincide with $\,p\sptilde\,$ as distributions.
Also, let ${\,w_{\,r}\,\in\,[\,0\,,\,1\,]\,}$; 
${\,r\in\,B^{p}\,}$, be a set of \emph{weights}   
such that
${\,\sum\nolimits_{\,r\in B^{p}} w_r\,=\,1\,}$. 

Then, every DSC field which in the scattering channel representation
on each component ${\,q\in\,B^{p}\,\subset\,\partial\zeta\,}$ equals 
${\,(\,z^{\,q}\,,\;z^{\,q\sptilde})\,}$ with
\vspace{-4pt}
\begin{equation}\centering\label{2.1}
z^{\,q\sptilde}\quad=\quad\sum\nolimits _{\,r\in B^{p}}\,w_r\;z^{\,r}
\vspace{-4pt}
\end{equation}
is named an (\,\emph{in $\zeta\,$ and with weights $w_{\,r}\,$}\,) 
\emph{coarse-grained field} (\,\emph{pertinent                 \linebreak
to $\,p\,$} -~or also \emph{to its nodal image $\,p\sptilde\,$\,}).
\end{udefinition}

The nodal state of every such field, coarse-grained in any cell, is hence
a convex superposition -~i.a. a \emph{weighted mean\,}~- of states
that represent the field on the cell surface.
Also, substituting a DSC field with any coarse-grained field
that (obviously) leaves the cell boundary states unchanged
is then called \emph{coarsening this field\;} (in any component). 

Should such an innocent looking procedure \emph{untie the Gordian knot 
of turbulence\/} ? Remember that our modest aim is repairing some defects
of discretization -~rather than resolving the fundamental questions of
turbulence in the way mathematical physics should do that.
What matters here is that in the DSC setup of CFD well-timed periodic 
coarse-graining of the flow efficiently prevents turbulent pile-up,
all with retaining the large eddies.
Note that the coarsening period should be taken as small as necessary
to ensure stability and precision, but large enough against the time
step (by one order of magnitude at least -~as a rule of thumb):
Clearly, coarsening inevitably also interferes with the flow dynamics.
However -~since the \emph{skeleton field\/} on the cell boundaries
is left unchanged~- the perturbation is minimal, and in fact negligible,
for suitable periods.
Through removing only the turbulent fluctuations of very short range 
(below a scale given by cell size) which tend to artificially pile-up,
coarsening regularises the flow to just such a mild degree that its essential
\emph{large} patterns are preserved.

%% file: sec3.tex
\section{Viscous Boussinesq-incompressible flow}\label{S:sec3}
The port and node distributions of a DSC algorithm can be
\emph{finite integrals}
-~as is the case for the TLM method, where finite path integrals over
electric and magnetic fields are evaluated in a discrete approximation
to Maxwell's integral equations [He4-6].
In the case at hand, they are simply \linebreak \emph{Dirac measures} that
pointwise evaluate the fields within the cells and on their surfaces,
and composites of Dirac measures which approximate the gradients of these 
fields; cf.~section~\eqref{S:sec4}. 

In the \emph{Oberbeck-Boussinesq\/} approximation \cite{Obb,Bss}
all material fluid properties are assumed as constant
-~\emph{except fluid density which only in the gravitational term
varies linearly with temperature}.
The energy equation of a \emph{Boussinesq-incompressible\/} fluid
of velocity $\,\vec{u}\,$ with thermal diffusivity $\,\alpha\,$, 
heat source(s) $\,q\,$,
and negligible viscous dissipation and heat capacity is then the
\emph{convection-diffusion} equation for the temperature $\,T\,$, 
cf.~\cite{GDN},
\begin{equation}\centering\label{3.1}
\frac{\partial\,T}{\partial\,t}\,
+\,div\,(\,T\,\vec{u}\,)\;=\;\alpha\,\Delta\,T\;+\;q\;.
\end{equation}
The Navier-Stokes \emph{momentum equations\/} for (\,the k-th component of\,)
a Newtonian fluid of dynamic viscosity $\,\eta\,$,
under pressure $\,p\,$,
and in a gravitational field of acceleration $\vec{g}$ 
take the form [e.g. Ptk]
\begin{equation}\centering\label{3.2}
\frac{\partial}{\partial\,t}\,(\,\varrho\,\vec{u}_{\,k}\,)\;
+\,div\,(\,\varrho\,\vec{u}_{\,k}\,\vec{u}\,)\;
+\;grad\,p_{\,k}\,
=\;\eta\,\Delta\,\vec{u}_{\,k}\;
+\;\varrho\;\vec{g}_{\,k}\;r.
\end{equation}
With
$\,\varrho_{\,\infty}=const\,$ and
$\,\varrho\,(\,T\,)=\varrho_{\,\infty}\,\beta\,
(\,T\,(\,t\,,\vec{x}\,)\,-\,T_{\,\infty}\,)\;$, wherein 
$\beta\,=\,\varrho^{-1}\,\partial\,\varrho\,/\,\partial\,T$,
this becomes in the OB-approximation
\begin{equation}\centering\label{3.3}
\frac{\partial\,\vec{u}_{\,k}}{\partial\,t}\,
+\,div\,(\,\vec{u}_{\,k}\,\vec{u}\,)\,
+\,\frac{grad\,p_{\,k}}{\varrho_{\infty}}\,
=\,\frac{\eta}{\varrho_{\infty}}\,\Delta\,\vec{u}_{\,k}\;
+\,\beta\,(\,T\,(\,t,\,\vec{x}\,)
- T_{\infty}\,)\,\vec{g}_{\,k}\,.
\end{equation}
Integrating these equations over cell $\zeta$ with boundary  
${\partial\zeta\;}$
and applying Gauss' theorem to integrals over 
$\,{\Delta}\,=\,{div\,grad\,}$ and 
$\,div\,(\,f\,\vec{u}\,)\,$
yields with time increment $\tau\,$ 
the following updating instructions for nodal $T$ and $\vec{u}_{\,k}$\,,
these quantities and $q\,$ averaged over the cell volume
${V_{\zeta}}$
\vspace{-4pt}
\begin{equation}\centering\label{3.4}
\begin{aligned}
T\,(\,t+\frac{\tau}{2}\,)\,
&=\;T\;+\\
&+\;\frac{\tau}{V_{\zeta}}\,
\int\nolimits_{\partial\,\zeta}(\,\alpha\,grad\,T\,
-\;T\,\vec{u}\,)\,\cdot\,dS
\end{aligned}
\vspace{-10pt}
\end{equation}
and
\vspace{-10pt}
\begin{equation}\centering\label{3.5}
\begin{aligned}
\vec{u}_{\,k}\,(\,t+\frac{\tau}{2}\,)\,
&:\,=\;\vec{u}_{\,k}\,
+\;\tau\,(\,\beta\,(\,T\,-\,T_{\infty}\,)\,\vec{g}_{\,k}\,
-\,\frac{grad\,p_{\,k}}{\varrho_{\infty}}\,\,)\;+\\
&+\;\frac{\tau}{V_{\zeta}}\;
\int\nolimits_{\partial\,\zeta}\,(\,\frac{\eta}{\varrho_{\infty}}\
grad\,\vec{u}_{\,k}\,
-\,\vec{u}_{\,k}\,\vec{u}\,)\,\cdot\,dS\;.
\end{aligned}
\end{equation}
In equations \eqref{3.4} and \eqref{3.5} the former updates
(at time $t-\tau/2\/$) of the \emph{nodal} quantities enter the right-hand
sides in the \emph{first line}, and the last former updates (at time $t\/$)
of the \emph{cell face} quantities enter the \emph{second line}\/.

The nodal values of $\,T\,$ and $\,\vec{u}_{\,k}$ are updated
at the reflection step 
-~after coarsening the velocity field, e.g. with weights
proportional to the pertinent cell face areas,
whenever the coarsening period is attained~-
while the port quantities that enter the surface integrals at the
right-hand sides are updated on the connection step of the iteration cycle.
The next section outlines how we can proceed with that 
in a non-orthogonal hexahedral mesh.

%% file: sec4.tex
\section{The non-orthogonal hexahedral cell}\label{S:sec4}
The physical interpretation of a DSC algorithm associates a
smoothly varying (\/e.g. in time and space $\,C^{\,\infty}$-\,)
scalar or vector field $\,Z\,$ to port and node states
$\,z^{\,p}\,$ and $\,z^{\,n}\,$ of a mesh cell system.

Let any hexahedral cell be given by its eight vertices.
Define then \emph{edge vectors} ${(_{\nu} e)_{\nu=0,...,11}}\/$,
\emph{node vectors\/} ${(_{\mu} b)_{\mu = 0,1,2}}$,
and \emph{face vectors\/} ${(_{\iota} f)_{\iota = 0,...,5}}$,\,
using the labelling scheme of figure~\ref{F:1}\,a
\vspace{-4pt}
\begin{equation}\label{4.1}\centering
\begin{split}
\begin{aligned}
_{\mu}b\;&:\,=\quad\frac{1}{4}
&& \!\!\sum\nolimits_{\nu = 0}^{3}\,_{_{(4\mu+\nu)}}e\,\,
&&\mu\,=\,0,1,2\\
\text{and}\quad
_{\iota}f\;&:\,=\;\,\frac{(-1)^{\,\iota}}{4}
&&\,(\,\,_{_{(8+2\iota)}}e\,
+_{_{(9+2(\iota+(-1)^{\iota}))}}e\,)
\,\,\land &&\\
& &&\;\;\land\,(\,_{_{(4+2\iota)}}e\, 
+_{_{(5+2\iota)}}e\,)\,\,&&\iota\,=\,0,...,5\,,
\end{aligned}
\end{split}
\vspace{-3pt}
\end{equation}
with all indices understood cyclic modulo 12\, and $\,\land\,$
denoting the cross product in $\mathbb{R}^3$.
\vspace{-25pt}
\input{fig23.tex}
At every cell face ${\iota\in\{0,...,5\}}\/$ and for any 
given $\tau\in\mathbb{R}_{+}\,$ the following time shifted finite
differences of $\,Z\,$ in directions ${ _{\mu} b }$ (\,$\mu = 0,1,2\,$)
form a vector valued function
\vspace{-4pt}
\begin{equation}\label{4.2}\centering
\begin{split}
_{\iota}\!{\nabla}^{B} Z_{\mu}\,(\,t\,)\;:\,=\;
\begin{cases}
\,2\,(-1)^{\iota}(\,Z^{n}\, _{\mid\,t-\tau/2} 
-\,_{\iota}Z^{p}\,_{\mid\,t}\,)\quad
&\text{if $\mu\,=\,[\iota / 2]$}\\
\,(\,\,_{2\mu+1} Z^{p}\,-\,_{2\mu} Z^{p}\,\,)
\,_{\mid\,t -\tau\,}\quad
&\text{if $\mu\,\neq\,[\iota / 2]$}\,
\end{cases}
\end{split}
\vspace{-4pt}
\end{equation}
($\,[\,x\,]$ denotes the \emph{integer part} of $\,x\in\mathbb{R}\,$).
The time increments are chosen conform with the updating conventions
of DSC schemes (as will be seen in a moment) and are consistent.
In fact, in the first order of the time increment ${\,\tau\,}$
and of the linear cell extension,
the vector ${\,_{\iota}\!{\nabla}^{B} Z\,}$
in the centre point of face ${\,\iota\,}$ approximates
the scalar products of the node vectors with the gradient
${\,\nabla Z\,}$.
More precisely, let for a fixed centre point on face ${\,\iota\,}$
and $\,\epsilon \in \mathbb{R}_{+}\,$ the \emph{$\epsilon$-scaled cell}
have edge vectors
$\,_{\iota}e\sptilde\,:\,=\,\epsilon\,\,_{\iota}e\,$. 
Let also $\,_{\iota}{\nabla}^{B\sptilde}Z_{\mu}\,$ denote
function \eqref{4.2} for the $\epsilon$-scaled cell (with node vectors 
$\,_{\mu}b\sptilde\,=\,\epsilon\,_{\mu}b\;$).
Then at the fixed point holds
\vspace{-4pt}
\begin{equation}\label{4.3}\centering
\begin{split}
<\,_{\mu}b\,,\,\text{grad($Z$)}\,>\,\,\,
=\,\,_{\mu}b\cdot\nabla Z\,\,
=\,\,\lim_{\epsilon\to 0}\,\,\lim_{\tau \to 0} \, \,
\frac{1}{\epsilon}\,_{\iota}\!{\nabla}^{B^{\sptilde}}Z_{\mu} \, ,
\end{split}
\vspace{-4pt}
\end{equation}
as immediately follows from the required $C^1$-smoothness of the field $Z$.

To recover the gradient
${{\nabla}Z\,}$ from \eqref{4.2}
in the same order of approximation,
observe that for every orthonormal
basis $\,{(_{\nu}u)_{\nu=0,...,m-1}}\,$ of
$\mathbb{R}^{m}\,\text{or}\,\,\mathbb{C}^{m}\,$, and for
any basis $\,{(_{\mu}b)_{\mu = 0,...,m-1}}\,$ with coordinate
matrix ${\beta_{\nu}^{\mu}}\,=\,{<\,_{\nu}u\,,\, _{\mu}b\,>}$,
the scalar products of every vector $\,a\,$ with $\,{_{\mu} b}\,$ equal
\vspace{-4pt}
\begin{equation}\label{4.4}\centering
\underbrace{<\,_{\mu}b\,,\,a\,>}_{\qquad=\,: 
\,\,{\alpha}_{\mu}^{B}}\,\,=\,\sum\nolimits_{\nu=0}^{m-1}\,
\underbrace{<\,_{\mu}b\,,\,_{\nu}u\,>}_{\,\,\,
({\bar{\beta}}_{\mu}^{\nu})\,=\,({\beta}_{\nu}^{\mu} )^{^{*}}}\,
\underbrace{<\,_{\nu}u\,,\,a\,>}_{\qquad=\,:\,\,{\alpha}_{\nu}}
\,\,=\,\bar{\beta}_{\mu}^{\nu}\,{\alpha}_{\nu}\;
\vspace{-4pt}
\end{equation}
(\,at the right-hand side, and henceforth, we observe
\textsc{Einstein}'s summation convention -~yet without summing
over indices that anywhere appear also as a left-hand subscript\,),
hence
\vspace{-6pt}
\begin{equation}\label{4.5}\centering
{\alpha}_{\nu}\,=\,
{\gamma}_{\nu}^{\mu}\alpha_{\mu}^{B}\,,
\qquad\text{with}\qquad({\gamma}_{\nu}^{\mu}) 
\,:\,=\,{({(\beta_{\nu}^{\mu})}^{*})}^{-1}\quad .
\vspace{-4pt}
\end{equation}
This applied to the node vector basis ${ _{\mu}b\,}$ and \eqref{4.3}
yields the approximate gradient of $\,Z\,$ at face $\iota$
\vspace{-6pt}
\begin{equation}\label{4.6}\centering
_{\iota}\!\nabla Z_{\nu}\quad
=\quad{\gamma}_{\nu}^{\mu}\,\,\,_{\iota}\!{\nabla}^{B}Z_{\mu}.
\vspace{-4pt}
\end{equation}
The scalar product of the gradient with face vector  
${_{\iota}f^{\nu}}\,=\,{<\,_{\iota}f,\,_{\nu}u>}\,$,\\
$\,\nu\in\{0,1,2\}\,$ is thus
\vspace{-6pt}
\begin{equation}\label{4.7}\centering
_{\iota}S\; 
=\;_{\iota}f\,\cdot\,_{\iota}\!\nabla Z\;
=\!\!\underbrace{_{\iota}f^{\nu}\;
{\gamma}_{\nu}^{\mu}}_{\qquad\;=\,:\,\,_{\iota} s^{\mu}}
\!\!_{\iota}\!{\nabla}^{B} Z_{\mu}\,\,
=\,\,_{\iota} s^{\mu}\,\,_{\iota}\!{\nabla}^{B} Z_{\mu}\;.
\vspace{-4pt}
\end{equation}
Continuity of the gradient at cell interfaces yields
linear updating equations for $Z^p$ on the two adjacent faces.  
In fact, for any two neighbouring cells
$\zeta$, $\chi$ with common face, labelled $\iota$
in cell $\zeta$ and $\kappa$ in $\chi\,$, continuity requires
\vspace{-6pt}
\begin{equation}\label{4.8}\centering
_{\iota}^{^{\zeta}}\!S\quad=\quad-\,\,_{\kappa}^{^{\chi}}\!S\,.
\vspace{-4pt}
\end{equation}
Substituting \eqref{4.7} for
$\,_{\iota}^{^{\zeta}}\!S\,$ and $\,_{\kappa}^{^{\chi}}\!S\,$
and observing the time shifts in \eqref{4.2}
\linebreak
provides the updating relations for 
$\,Z^{\,p}\,$ at the cell interfaces.
To make these explicit, we first introduce the quantities
$\,{_{\iota}z_{\mu}^{\,p,\,n}}\,$, ($\,\iota\,=\,0,...,5\,$;
$\mu\,=\,0,1,2\,$)
\vspace{-6pt}
\begin{equation}\label{4.10}\centering
\begin{split}
_{\iota}z_{\mu}^{n}\,(\,t\,)\quad :\,=\quad
\begin{cases}
\,\,2\,(-1)^{\iota}\,\,Z^{\,n}\,_{\mid\,t}\qquad
&\text{if $\mu\,=\,[\iota /2]$}\,\,\\
\,\,(\,_{2\mu +1} Z^{\,p} 
-\,_{2\mu} Z^{\,p}\,)_{\mid\,t-\tau/2}\qquad
&\text{else}\,
\end{cases}\;,
\end{split}
\vspace{-4pt}
\end{equation}
which in virtue of \eqref{1.1} yields
$\;{_{\iota}z_{\mu}^{\,p}}\,=\,{(\,p\,,\,Z\,)}\,
=\,{(\,p\sptilde,\,Z\circ\,_{\iota}{\sigma}^{-1}\,)}\,
=\\=\,{_{\iota}z_{\mu}^{\,n}\,\mid{Z\,\circ\,_{\iota}{\sigma}^{-1}}}\,$,
where $\,_{\iota}{\sigma}\,:\,n\,\mapsto\,p\,$
denotes the nodal shift pertinent to face $\,\iota\,$.\;
In particular
\vspace{-4pt}
\begin{equation}\label{4.11}
_{\iota}z_{[\iota/2]}^{\,p}\,(\,t\,)\quad=\quad
\,\,2\,(-1)^{\iota}\,\,_{\iota}Z^{\,p}\,_{\mid\,t}\;,
\end{equation}
which together with \eqref{4.10} is consistent for 
$\,\mu\,\neq\,[\iota/2]\,$ with
\vspace{-4pt}
\begin{equation}\label{4.12}
_{\iota}z_{\mu}^{\,n}\,(\,t+\tau/2\,)\quad
=\quad -\;\frac{1}{2}\,(\,_{\,2\mu+1}z_{\mu}^{\,p}\,
+\, _{\,2\mu}z_{\mu}^{\,p}\,)\,(\,t\,)\;.
\vspace{-4pt}
\end{equation}
From (\,\ref{4.2}, \ref{4.7}, \ref{4.10}, \ref{4.11}\,) follows that
\vspace{-4pt}
\begin{equation}\label{4.13}
\begin{split}
_{\iota}S\,_{\mid\,t+\tau}\quad
&=\quad\,_{\iota}s^{\mu}\,(\,_{\iota}z_{\mu}^{\,n}\,_{\mid\,t+\tau/2}\,
-\,2\,{(-1)}^{\iota} {\delta}_{\mu}^{[\iota/2]}\,\,
_{\iota}Z^{\,p}\,_{\mid\,t+\tau}\,) \\
&=\quad _{\iota}s^{\mu}\,(\,_{\iota}z_{\mu}^{n}\,_{\mid\,t+\tau/2} \,
-\,{\delta}_{\mu}^{[\iota/2]}\,\, 
_{\iota}z_{\mu}^{\,p}\,_{\mid\,t+\tau}\,)
\; .
\end{split}
\vspace{-4pt}
\end{equation}
Continuity of $\,Z\;$, i.e.
$\,_{\iota}^{^{\zeta}}Z\,^{p}\,=\,_{\kappa}^{^{\chi}}Z\,^{p}\,$,\, 
with (\ref{4.8}, \ref{4.10}) then implies
\vspace{-6pt}
\begin{equation}\label{4.14}
_{\iota}^{^{\zeta}}z\,_{[\iota/2]}^{p}\,(\,t+\tau\,)\,=\,
\,\frac{\,_{\iota}^{^{\zeta}}s\,^{\mu}\,\, 
_{\iota}^{^{\zeta}}z\,_{\mu}^{n}\,(\,t+\tau/2\,)\,
+\,_{\kappa}^{^{\chi}}s\,^{\nu}\,\,\,
_{\kappa}^{^{\chi}}z\,_{\nu}^{n}\,(\,t+\tau/2\,)} 
{_{\iota}^{^{\zeta}}s\,^{[\iota/2]}
+\,( -1 )^{\iota +\kappa}\;_{\kappa}^{^{\chi}}s\,^{[\kappa/2]}}\;.
\vspace{-4pt}
\end{equation}
For completeness we agree upon setting
$\,_{\iota}^{^{\zeta}}z\,_{\mu}^{p}\,(\,t+\tau\,)\,
:\,=\,_{\iota}^{^{\zeta}}z\,_{\mu}^{n}\,(\,t+\tau/2\,)\,$
for $\,\mu\,\neq\,[\iota/2]\,$
(although this contains a slight inconsistency
in that continuity might be infringed; this can easily be remedied
by taking the arithmetic mean of the two adjacent values).
 -~In fact, our agreement doesn't do harm, since any discontinuity
disappears with mesh refinement. 

We dispose, hence, of a complete set of recurrence relations for
$\,{z^{\,p}}\,$ (\,given $\,{z^{\,n}}\,$ by the former reflection step)
which at the same time determine
the field components on face~$\iota$ and their gradients
\vspace{-6pt}
\begin{equation}\label{4.15}\centering
_{\iota}\!\nabla Z_{\nu}\quad
=\quad{\gamma}_{\nu}^{\mu}\; _{\iota}z_{\mu}^{\,p}.
\vspace{-4pt}
\end{equation}
Essentially this constitutes the connection step of the algorithm.
\newline
Nodal gradients are similarly (and even more simply) derived, using
\vspace{-6pt}
\begin{equation}\label{4.16}\centering\notag
{\nabla}^{B} Z_{\mu}^{\,n}\,(\,t\,+\frac{\tau}{2}\,)\;\;:\,
=\;\;(\,_{2\mu+1} Z^{\,p}\,-\,_{2\mu} Z^{\,p}\,)\,(\,t\,)\;;
\quad\mu\,=\,0,\,1,\,2\;
\vspace{-4pt}
\end{equation}
in the place of \eqref{4.2}, and then again \eqref{4.6}.
With the node and cell-boundary values and gradients of $\,T\,$ and
$\,\vec{u}\,$ the nodal updating relations for the latter are
immediately extracted from equations~(\,\ref{3.4}, \ref{3.5}\,)
in section~\ref{S:sec3}\,. For equation \eqref{3.4} this is
essentially (up to the convective term) carried out in
\cite{He1}, section~5, and the procedure remains straightforward
in the case at hand.
Note that a well-timed LES coarsening routine, cf. section~\ref{S:sec2},
should be periodically carried out before the nodal step of the iteration
cycle in order to obviate instablities called forth by the energy 
cascade~\cite{Po}. 

The updating relations thus obtained are explicit and consistent
with near-field interaction (\,only adjacent quantities enter\,).
So, they can optionally be transformed into scattering relations
for incident and reflected quantities \eqref{1.6} along the
guidelines of section~\ref{S:sec1} -~with established advantages
for stability estimates \cite{He2}.

%% file: fig23.tex
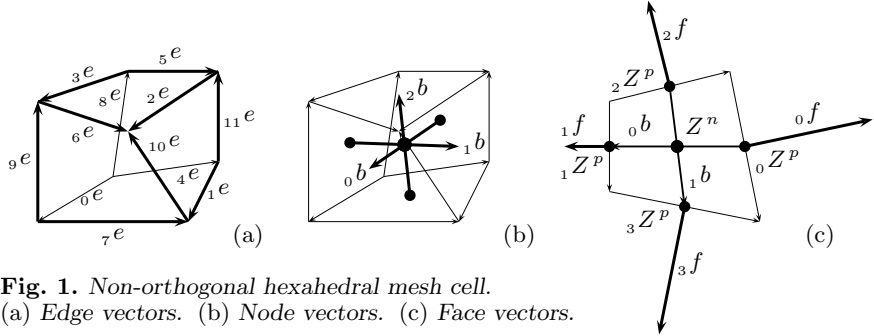
\begin{figure}[!h]\centering
\setlength{\unitlength}{1.cm}
\begin{pspicture}(-.9,-.5)(20,3.5)\centering
\psset{xunit=.4cm,yunit=.4cm}
\psline[linewidth=0.1mm]{->}(1.5,1.5)(-1.0,0.0)
\psline[linewidth=0.1mm]{->}(1.5,1.5)(5.0,2.0)
\psline[linewidth=0.1mm]{->}(1.5,1.5)(2.0,5.0)
\psline[linewidth=0.4mm]{->}(-1.0,0.0)(-1.0,4.0)
\psline[linewidth=0.4mm]{->}(-1.0,0.0)(4.0,0.0)
\psline[linewidth=0.4mm]{->}(5.0,2.0)(5.0,5.0)
\psline[linewidth=0.4mm]{->}(5.0,2.0)(4.0,0.0)
\psline[linewidth=0.4mm]{->}(4.0,0.0)(2.0,3.0)
\psline[linewidth=0.4mm]{->}(2.0,5.0)(-1.0,4.0)
\psline[linewidth=0.4mm]{->}(2.0,5.0)(5.0,5.0)
\psline[linewidth=0.4mm]{->}(-1.0,4.0)(2.0,3.0)
\psline[linewidth=0.4mm]{->}(5.0,5.0)(2.0,3.0)
\rput(0.8,0.8){$_{_{0}}e$}
\rput(5.0,1.0){$_{_{1}}e$}
\rput(3.0,4.2){$_{_{2}}e$}
\rput(0.5,4.9){$_{_{3}}e$}
\rput(4.0,1.5){$_{_{4}}e$}
\rput(3.4,5.5){$_{_{5}}e$}
\rput(0.5,2.9){$_{_{6}}e$}
\rput(1.5,-.5){$_{_{7}}e$}
\rput(1.4,4.1){$_{_{8}}e$}
\rput(-1.6,2.0){$_{_{9}}e$}
\rput(3.2,2.6){$_{_{10}}e$}
\rput(5.7,3.4){$_{_{11}}e$}
\rput(6.0,-0.5){\small{\textnormal{(a)}}}
\psline[linewidth=0.1mm]{->}(10.5,1.5)(08.0,0.0)
\psline[linewidth=0.1mm]{->}(10.5,1.5)(14.0,2.0)
\psline[linewidth=0.1mm]{->}(10.5,1.5)(11.0,5.0)
\psline[linewidth=0.1mm]{->}(08.0,0.0)(08.0,4.0)
\psline[linewidth=0.1mm]{->}(08.0,0.0)(13.0,0.0)
\psline[linewidth=0.1mm]{->}(14.0,2.0)(14.0,5.0)
\psline[linewidth=0.1mm]{->}(14.0,2.0)(13.0,0.0)
\psline[linewidth=0.1mm]{->}(13.0,0.0)(11.0,3.0)
\psline[linewidth=0.1mm]{->}(11.0,5.0)(08.0,4.0)
\psline[linewidth=0.1mm]{->}(11.0,5.0)(14.0,5.0)
\psline[linewidth=0.1mm]{->}(08.0,4.0)(11.0,3.0)
\psline[linewidth=0.1mm]{->}(14.0,5.0)(11.0,3.0)
\psline[showpoints=true,linewidth=0.4mm]{->}(12.375,3.375)(10.0,1.75)
\psline[showpoints=true,linewidth=0.4mm]{->}(09.375,2.625)(13.0,2.50)
\psline[showpoints=true,linewidth=0.4mm]{->}(11.37,.875)(11.0,4.25)
\psline[showpoints=true,
linewidth=0.6mm]{-}(11.1875,2.5625)(11.1875,2.5625) 
\rput(09.55,1.55){$_{_{0}}b$}
\rput(13.5,2.60){$_{_{1}}b$}
\rput(11.6,4.40){$_{_{2}}b$}
\rput(15,-0.5){\small{\textnormal{(b)}}}
\psline[linewidth=0.1mm]{->}(18.0,4.0)(18.0,1.0)
\psline[linewidth=0.1mm]{->}(18.0,1.0)(23.0,0.0)
\psline[linewidth=0.1mm]{->}(18.0,4.0)(22.0,5.0)
\psline[linewidth=0.1mm]{->}(22.0,5.0)(23.0,0.0)
\psline[showpoints=true,linewidth=0.5mm]{-}(20.25,2.5)(20.25,2.5) 
\psline[linewidth=0.25mm]{->}(22.5,2.5)(18.0,2.5)
\psline[linewidth=0.25mm]{->}(20.0,4.5)(20.5,0.5)
\psline[showpoints=true, linewidth=0.4mm]{->}(22.5,2.5)(26.751,3.400)
\psline[showpoints=true,linewidth=0.4mm]{->}(18.0,2.5)(16.5,2.5)
\psline[showpoints=true,linewidth=0.4mm]{->}(20.0,4.5)(19.287,7.352)
\psline[showpoints=true,linewidth=0.4mm]{->}(20.5,0.5)(19.650,-3.751)
\rput(19.0,3.2){$_{_{0}}b$}
\rput(21.0,1.5){$_{_{1}}b$}
\rput(24.6,3.6){$_{_{0}}f$}
\rput(16.8,3.2){$_{_{1}}f$}
\rput(20.2,6.4){$_{_{2}}f$}
\rput(20.7,-1.4){$_{_{3}}f$}
\rput(23.6,2.0){$_{_{0}}Z^{\,p}$}
\rput(17.1,1.8){$_{_{1}}Z^{\,p}$}
\rput(18.8,4.7){$_{_{2}}Z^{\,p}$}
\rput(19.3,0.0){$_{_{3}}Z^{\,p}$}
\rput(21.1,3.2){$Z^{\,n}$}
\rput(25,-0.5){\small{\textnormal{(c)}}}
\end{pspicture}
\caption{\textsl{Non-orthogonal hexahedral mesh cell. \newline
\textnormal{(a)} Edge vectors.\;
\textnormal{(b)} Node vectors.\;
\textnormal{(c)} Face vectors.}\hfill}\label{F:1}
\end{figure}

%% file: sec5.tex
\section{Pressure}\label{S:sec5}
Unlike compressible flow -~wherein a thermodynamic \emph{state equation}
relates pressure to density and temperature~-
Boussinesq-incompressible flow is conserved by pressure acting
like a potential against violations of the local flow balance.
Hence, the pressure is directly coupled to the flow divergence by
\emph{Poisson's equation}:
\vspace{-0.15cm}
\begin{equation}\centering\label{5.1}
\Delta p\,=\,(\,\varrho_{\infty}/\tau\,)\,div\,\vec{u}\,.
\vspace{-0.10cm}
\end{equation}
The solution $p$ for fluid velocities updated in the connection step
provides a pressure gradient, which in the next reflection step 
-~via equations (\ref{3.3},\ref{3.5})~- repairs sporadic violations
of the flow balance.
Systematic test computations have shown that additional 
\emph{divergence clearing}, as proposed by the author in earlier
papers, e.g.~\cite{He3}, is not necessary (nor even profitable)
-~and may be rather on the debit of computational performance and stability.

In integral form, after using Gauss' Theorem, equation \eqref{5.1} becomes
\vspace{-0.15cm}
\begin{equation}\centering\label{5.2}
\quad\tau\,
\int\nolimits_{\,\partial\zeta}\, grad\,p\,\cdot\,dF\quad
=\quad\varrho_{\infty}\,
\int\nolimits_{\,\partial\zeta}\,\vec{u}\,\cdot\,dF\;,
\vspace{-0.10cm}
\end{equation}
which can be solved by a Gauss-Seidel routine or by successive
overrelaxation, carried out after updating the face velocities
in the connection step.

In each iteration, firstly the discrete cell boundary integral 
${\,I_{\partial\zeta}}\,: =\,{\int _{\partial\zeta}\vec{u}\cdot dF\,}$
is computed and explicitely -~by solving a linear equation in 
$\,I_{\partial\zeta}\,$ and the cellular pressures~- the unique $\,p^{n}\,$
that exactly compensates $\,I_{\partial\zeta}\,$ 
so that \eqref{5.2} is satisfied in every cell. In a second run
over the mesh, continuity of the pressure gradient is re-established
at the cell interfaces by updating the face pressures along the
lines of the last section (viz. taking $\,Z\,$ as $\,p\,$ there). 
The loop is reiterated until (usually after a few iterations) equations
\eqref{5.2} hold with sufficient precision.

%% file: ccl.tex
\section{Conclusion and Completions}\label{S:ccl}
The DSC approach to computational fluid dynamics allows for
simulation of turbulent flow in a plain and natural way.
Cellular coarse-graining efficiently prevents turbulent pile up
and preserves the large eddies. 

Since the energy cascade starts long before transition to turbulence
(actual\-ly it is effective over the entire energy spectrum~\cite{Po}, 
\cite[pp.~72~ff.]{BIL}\,), a well-timed coarsening routine should 
escort every DSC fluid flow algo\-rithm. 
Periodic coarse-graining significantly improves algorithm stability,
even in the laminar regime -~this without overly altering the flow dynamics,
if the coarsening period is chosen as outlined in section~\ref{S:sec2}.  
Handy criteria for a good choice -~in relation to mesh refinement,
e.g.~- should be subject to further study. 

Details of implementation and special applications are clearly 
beyond the scope of the present study. 
However, we leave off with some examples for the purpose of illustration. 

The graphics displayed have been computed with author's test program
\textsc{DANSE}. The latter combines a TLM Maxwell field solver with a DSC
heat transfer and fluid flow algorithm written in the lines of this paper.

In particular, the computations for coaxial RF power transmission line RL100-230
(motivated by an ion cyclotron resonance heating experiment in plasma physics)
are in excellent keeping with reference computations (\textsc{Fluent}) 
and empirical data:
For a line with inner conductor made of copper and air dielectric at atmospheric 
pressure the heating process has been simulated in horizontal position
from standby to steady state operation with 160 KW transmitted CW power
at 100 MHz frequency, with outer conductor cooled at 40 degrees Celsius
temperature. Simultaneously, a Maxwell field TLM algorithm run
in the same mesh provided the skin effect heat sources.
\input{edge.tex}
\input{karman.tex}
\input{RL100-230_flow.tex}
\input{RL100-230_temp.tex}
\input{RL100-230_Tprofile.tex}

%% file: edge.tex
\begin{figure}[!h]\centering
\setlength{\unitlength}{1.0cm}
\begin{pspicture}(0.7,1.1)(12.9,3.5)\centering
\psset{xunit=0.9cm,yunit=0.9cm}
\rput(9.95,4.90){\rotateright{\includegraphics[scale=0.80,clip=0]{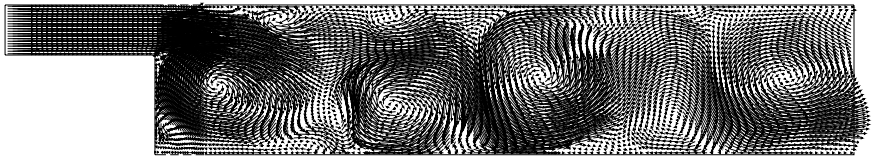}}}
\end{pspicture}
\caption{\textsl{Jet over backward facing step\qquad\qquad\qquad\qquad\qquad\newline
[ highly turbulent regime ].
\qquad\qquad\qquad\hfill}}
\label{F:2}
\end{figure}

%% file: karman.tex
\begin{figure}[!h]\centering
\setlength{\unitlength}{1.0cm}
\begin{pspicture}(0.7,1.2)(12.9,3.1)\centering
\psset{xunit=0.9cm,yunit=0.9cm}
\rput(9.95,2.50){\rotateright{\includegraphics[scale=0.77,
clip=0]{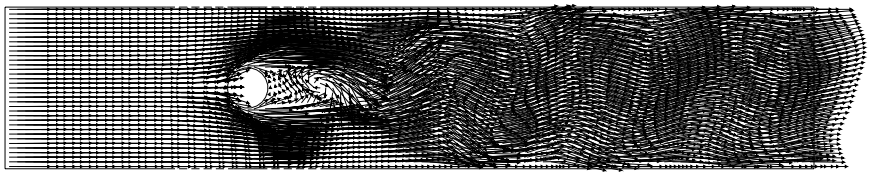}}}
\end{pspicture}
\caption{\textsl{K\'arm\'an vortex street behind cylinder\newline
[ snapshot of oscillating flow ].
\qquad\qquad\qquad\hfill}}
\label{F:3}
\end{figure}

%% file: RL100-230_flow.tex
\begin{figure}[!h]\centering
\setlength{\unitlength}{1.cm}
\begin{pspicture}(0.0,0.0)(12.0,8.7)\centering
\psset{xunit=0.9cm,yunit=0.9cm}
\rput(6.10,4.50){\rotateright{\includegraphics[scale=0.60,clip=0]{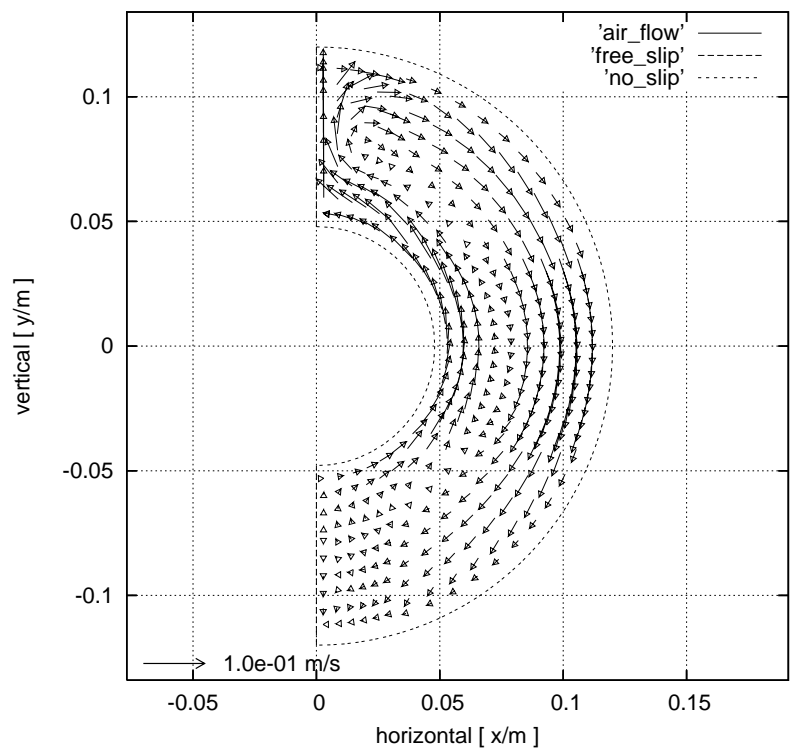}}}
\end{pspicture}
\vspace{25pt}
\caption{\textsl{Air convection in horizontal coaxial line RL100-230\;;
\newline
transverse section
[ reference arrow: 0.1$\,m s^{-1}\,$;
transmitted power 160 kW CW, frequency 100 MHz; air dielectric;
inner conductor copper, outer conductor cooled at 40 degrees Celsius ].
\qquad\qquad\hfill}}
\label{F:4}
\end{figure}

%% file: RL100-230_temp.tex
\begin{figure}[!h]\centering
\setlength{\unitlength}{1.cm}
\begin{pspicture}(0.0,0.0)(12.0,7.0)\centering
\psset{xunit=0.9cm,yunit=0.9cm}
\rput(6.50,3.67){\rotateright{\includegraphics[scale=0.47,clip=0]{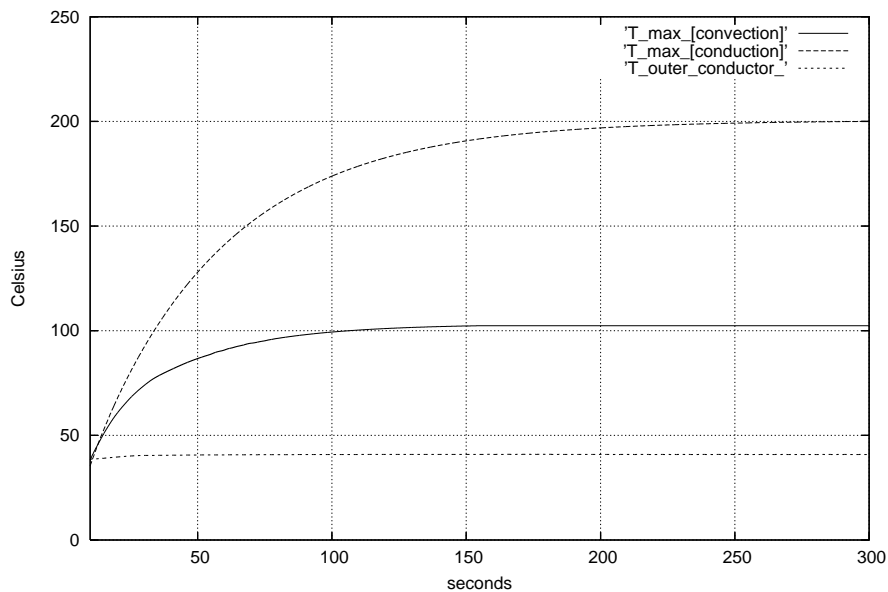}}}
\end{pspicture}
\vspace{15pt}
\caption{\textsl{Heating process in coaxial line RL100-230\;;\newline
temperature maximum vs. time from power-on,\newline
convection compared to -~fictitious~- pure conduction \newline
[ transmitted power 160 kW CW, frequency 100 MHz,\newline
outer conductor cooled at 40 degrees Celsius ].
\qquad\qquad\qquad\hfill}}
\label{F:5}
\end{figure}

%% file: RL100-230_Tprofile.tex
\begin{figure}[!h]\centering
\setlength{\unitlength}{1.cm}
\begin{pspicture}(0.0,0.0)(12.0,9.2)\centering
\psset{xunit=0.9cm,yunit=0.9cm}
\rput(5.00,4.50){\rotateright{\includegraphics[scale=0.60,clip=0]{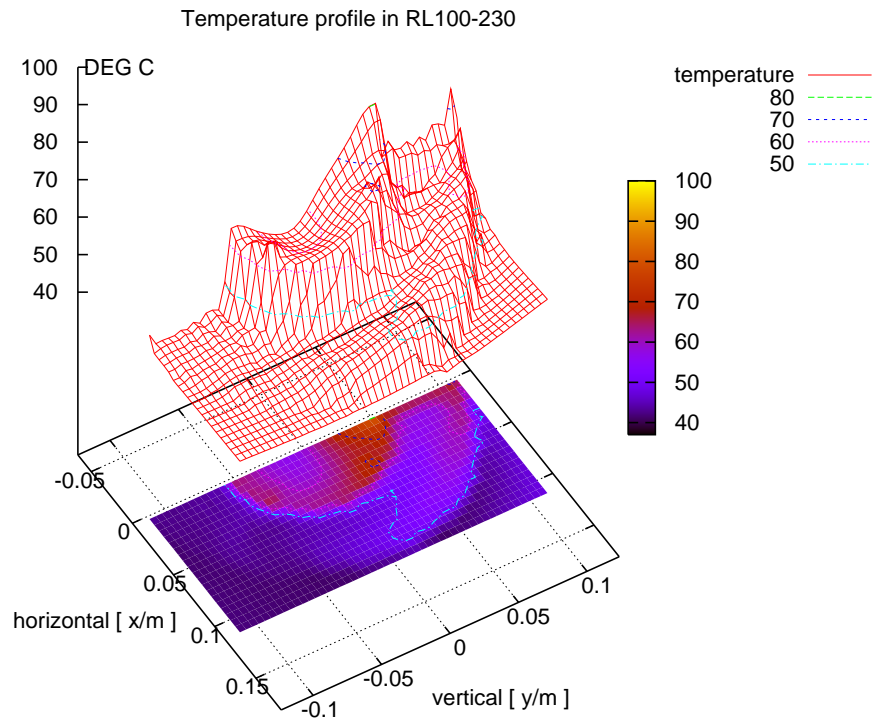}}}
\end{pspicture}
\vspace{10pt}
\caption{\textsl{Cross sectional temperature profile in coaxial line RL100-230
computed at steady state\;
[ transmitted power 160 kW CW, frequency 100 MHz; air dielectric;
inner conductor copper, outer conductor cooled at 40 degrees Celsius ].
\qquad\qquad\hfill}}
\label{F:6}
\end{figure}